\documentclass[10pt,twoside,english,american]{article}
\usepackage[T2A]{fontenc}
\usepackage[utf8]{inputenc}
\usepackage[a4paper]{geometry}
\usepackage{color}
\usepackage{mathrsfs}
\usepackage{amsmath}
\usepackage{amssymb}

\makeatletter
\usepackage[english,ukrainian]{babel}
\usepackage[active]{srcltx}
\usepackage{ifthen}
\usepackage{cite}
\usepackage{mathrsfs}
\usepackage{hyperref}
\usepackage{color}
\usepackage{amsthm}
\usepackage{cmap} %for selection and copying in pdf file

\newtheorem{theorem}{Theorem}
\newtheorem*{theoremA}{Theorem A}
\newtheorem*{theoremB}{Theorem B}
\newtheorem{lemma}{Lemma}
\newtheorem{assertion}{Assertion}
\newtheorem{nasl}{Corollary}
\newtheorem{properties}{Properties}
\newtheorem{definition}{Definition}
\newtheorem{denotation}{Notation}
\theoremstyle{remark}
\newtheorem{rmk}{Remark}
\newtheorem{example}{{\bf\em Example}}

\newcommand{\BeginDef}{\begin{definition}}
\newcommand{\EndDef}  {\end{definition}}
\newcommand{\BeginDenot}{\begin{denotation}}
\newcommand{\EndDenot}{\end{denotation}}
\newcommand{\BeginLem}{\begin{lemma}}
\newcommand{\EndLem}  {\end{lemma}}
\newcommand{\BeginThm}{\begin{theorem}}
\newcommand{\EndThm}  {\end{theorem}}
\newcommand{\BeginAs}{\begin{assertion}}
\newcommand{\EndAs}{\end{assertion}}
\newcommand{\BeginNasl}{\begin{nasl}}
\newcommand{\EndNasl}  {\end{nasl}}
\newcommand{\BeginProp}{\begin{properties}}
\newcommand{\EndProp}{\end{properties}}
\newcommand{\BeginRmk}{\begin{rmk}}
\newcommand{\EndRmk}  {\end{rmk}}
\newcommand{\BeginProof}{\begin{proof}}
\newcommand{\EndProof}{\end{proof}}
\newcommand{\BeginEx}{\begin{example}}
\newcommand{\EndEx}{\end{example}}

\newcommand{\doctitle}{}
\newcommand{\docinfoOdd}{}
\newcommand{\docinfoEven}{}
\newcommand{\docauthor}{}
\newcommand{\DOI}{}
\newcommand{\Preprint}{true}
\date{\empty}
\renewcommand{\@oddhead}{\docinfoOdd}
\renewcommand{\@evenhead}{\docinfoEven}

\makeatother

\usepackage{babel}
\begin{document}
\global\long\def\mathcircumflex{\mbox{\^{}}} 

\global\long\def\N{\mathbb{N}}
\global\long\def\Z{\mathbb{Z}}
\global\long\def\Q{\mathbb{Q}}
\global\long\def\R{\mathbb{R}}
\global\long\def\RRR#1{\widehat{\R^{#1}}}
\global\long\def\CC{\mathbb{C}}
\global\long\def\card{\mathbf{card}}
\global\long\def\diag{\textsf{\textbf{diag}}}
\newcommand{\nlm}{\nolimits}\global\long\def\Dm{\mathfrak{D}}
\global\long\def\Rg{\mathfrak{R}}
\global\long\def\T{\mathbf{T}}
\global\long\def\TT{\mathbb{T}}
\global\long\def\MM{\mathcal{M}}
\global\long\def\arc#1{#1^{[-1]}}
\global\long\def\fff{\mathop{\leftarrow}}
\global\long\def\nff{\mathop{\:\not\not\negmedspace\negmedspace\fff}}
\global\long\def\FFF{\mathop{\overleftarrow{\qquad\;}}}
\global\long\def\AAA{\mathcal{A}}
\global\long\def\BBB{\mathcal{B}}
\global\long\def\Bs{\mathfrak{Bs}}
\global\long\def\Ind#1{\mathcal{I}nd\left(#1\right)}
\global\long\def\ind#1{\mathtt{ind}\left(#1\right)}
\global\long\def\Lk#1{\mathcal{L}k\left(#1\right)}
\global\long\def\cL{\mathcal{L}}
\global\long\def\cQ{\mathcal{Q}}
\global\long\def\Tm{\mathbf{Tm}}
\global\long\def\TM{\mathbb{T}\mathbf{m}}
\global\long\def\BS{\mathbb{B}\mathfrak{s}}
\global\long\def\bs#1{\mathsf{bs}\left(#1\right)}
\global\long\def\BsB{\Bs(\BBB)}
\global\long\def\BSB{\BS(\BBB)}
\global\long\def\TmB{\Tm(\BBB)}
\global\long\def\TMB{\TM(\BBB)}
\global\long\def\tm#1{\mathsf{tm}\left(#1\right)}
\global\long\def\BSBB#1{\BS\left(\BBB_{#1}\right)}
\global\long\def\bbL{\mathbb{L}}
\global\long\def\bbK{\mathbb{K}}
\global\long\def\Ll{\bbL l}
\global\long\def\Ld{\bbL d}
\global\long\def\Lc{\bbL c}
\global\long\def\Lg{\bbL g}
\global\long\def\At{\mathcal{A}t}
\global\long\def\cU{\mathcal{U}}
\global\long\def\sU{\mathscr{U}}
\global\long\def\fU{\mathfrak{U}}
\global\long\def\sT{\mathscr{T}}
\global\long\def\vcB{\overleftarrow{\BBB}}
\global\long\def\vfU{\overleftarrow{\fU}}
\global\long\def\cZ{\mathcal{Z}}
\global\long\def\Zim#1{\cZ\mathsf{im}\left(#1\right)}
\global\long\def\bbU{\mathbb{U}}
\global\long\def\IndZ{\Ind{\cZ}}
\global\long\def\LkZ{\Lk{\cZ}}
\global\long\def\un#1#2{\left\langle #2\fff#1\right\rangle }
\global\long\def\unn#1#2#3{\un{#1}{!\:#2}#3}
\global\long\def\uni#1#2{\left\langle #1\rightarrow#2\right\rangle }
\global\long\def\ol#1{#1\mathcircumflex}
\global\long\def\oll{\ol{\lf}}
\global\long\def\w{\mathrm{w}}
\global\long\def\vv{\mathrm{v}}
\global\long\def\wx{\widetilde{\w}}
\global\long\def\x{\mathrm{x}}
\global\long\def\y{\mathrm{y}}
\global\long\def\fX{\mathfrak{X}}
\global\long\def\cX{\mathcal{X}}
\global\long\def\cY{\mathcal{Y}}
\global\long\def\lsL{\bbL}
\global\long\def\spX{\fX}
\global\long\def\bX{\mathbf{X}}
\global\long\def\Lp#1#2{\mathfrak{Lp}\left(#1,#2\right)}
\global\long\def\kpQ{\mathfrak{Q}}
\global\long\def\tpT{\mathcal{T}}
\global\long\def\rfrm#1{\mathfrak{#1}}
\global\long\def\lf{\rfrm l}
\global\long\def\mf{\rfrm m}
\global\long\def\pf{\rfrm p}
\global\long\def\kf{\rfrm k}
\global\long\def\lk#1#2{\mathbf{lk}_{#1}\left(#2\right)}
\global\long\def\fP{\mathfrak{P}}
\global\long\def\fS{\mathfrak{S}}
\global\long\def\bbS{\mathbb{S}}
\global\long\def\Zimm#1{\cZ\mathsf{im}\left[#1\right]}
\global\long\def\ES#1#2{#1^{\left\{  #2\right\}  }}
\global\long\def\es#1#2{#1_{\left\{  #2\right\}  }}
\global\long\def\ESl#1{\ES{#1}{\lf}}
\global\long\def\esl#1{\es{#1}{\lf}}
\global\long\def\ESm#1{\ES{#1}{\mf}}
\global\long\def\esm#1{\es{#1}{\mf}}
\global\long\def\omh{\hat{\omega}}
\global\long\def\het{\hat{\eta}}
\global\long\def\MS#1{\widehat{\mathbf{#1}}}
\global\long\def\MSA{\MS A}
\global\long\def\lsh#1{\mathscr{#1}}
\global\long\def\lshA{\lsh A}
\global\long\def\po#1{\mathsf{po}\left(#1\right)}
\global\long\def\ki#1{\mathsf{ki}\left(#1\right)}
\global\long\def\HDNN{\Uparrow}
\global\long\def\HDP{\HDNN^{+}}
\global\long\def\HDNP{\Downarrow}
\global\long\def\HDN{\HDNP^{-}}
\global\long\def\HDDf{\Updownarrow^{\pm}}
\global\long\def\Tpi{\mathfrak{Tpi}}
\global\long\def\esb{\underrightarrow{\subset}}
\global\long\def\esp{\underleftarrow{\supset}}
\global\long\def\Esb{\underrightarrow{\sqsubset}}
\global\long\def\Esp{\underleftarrow{\sqsupset}}
\global\long\def\Zk{\mathbf{Zk}}
\global\long\def\ZkQ{\Zk(\kpQ)}
\global\long\def\Tp{\mathcal{T}p}
\global\long\def\TpQ{\Tp(\kpQ)}
\global\long\def\Ls{\mathbb{L}s}
\global\long\def\LsQ{\Ls(\kpQ)}
\global\long\def\Ps{\mathfrak{Ps}}
\global\long\def\PsQ{\Ps(\kpQ)}
\global\long\def\ds{\mathbf{di}}
\global\long\def\dsQ{\ds_{\kpQ}}
\global\long\def\NQ#1{\left\Vert #1\right\Vert _{\kpQ}}
\global\long\def\ipQ#1#2{\left(#1,#2\right)_{\kpQ}}
\global\long\def\cG{\mathcal{G}}
\global\long\def\cK{\mathcal{K}}
\global\long\def\cG{\mathcal{G}}
\global\long\def\cK{\mathcal{K}}
\global\long\def\cC{\mathcal{C}}
\global\long\def\fC{\mathfrak{C}}
\global\long\def\bkK{\cK^{\mathfrak{b}}}
\global\long\def\bkC{\fC^{\mathfrak{b}}}
\global\long\def\cF{\mathcal{F}}
\global\long\def\cF{\mathcal{F}}
\global\long\def\obr#1#2{#1\upharpoonright#2}
\global\long\def\Cobr#1{\obr{\fC}{#1}}
\global\long\def\Col{\Cobr{\lf}}
\global\long\def\BE{\mathsf{BE}}
\global\long\def\BBE{\mathbb{BE}}
\global\long\def\BG{\mathsf{BG}}
\global\long\def\unc#1#2{\left[#2\fff#1\right]}
\global\long\def\LkF{\Lk{\cF}}
\global\long\def\tr#1#2{\mathfrak{trj}_{#1}\left[#2\right]}
\global\long\def\TR{\mathbb{T}\mathbf{rj}}
\global\long\def\TRd{\overline{\TR}}
\global\long\def\kr{\mathfrak{q}}
\global\long\def\Mk{\mathbb{M}k}
\global\long\def\sPk{\mathbf{Q}}
\global\long\def\pk#1{\sPk^{\left\langle #1\right\rangle }}
\global\long\def\PKK#1#2{\sPk^{\un{#1}{#2}}}
\global\long\def\gteq{\mathop{\rightleftarrows}}
\global\long\def\ngteq{\mathop{\not\rightleftarrows}}
\global\long\def\I{\mathbb{I}}
\global\long\def\Kim#1{\mathfrak{Kim}\left(#1\right)}
\global\long\def\Kimm#1{\mathfrak{Kim}\left[#1\right]}
\global\long\def\Ku#1{\mathfrak{Ku}\left[#1\right]}
\global\long\def\knu#1{\mathfrak{Ku}\left(#1\right)}
\global\long\def\LkC{\Lk{\fC}}
\global\long\def\bbUx{\widehat{\bbU}}
\global\long\def\pkeq{\left[\equiv\right]}
\global\long\def\frl#1#2{#1\downharpoonright_{#2}}
\global\long\def\visb{\,\underrightarrow{<}\,}
\global\long\def\visp{\,\underleftarrow{>}\,}
 \global\long\def\vie#1#2#3{#1^{\left\{  #2\rightarrow#3\right\}  }}
\global\long\def\Ha{\mathfrak{H}}
\global\long\def\MHa{\mathcal{M}\left(\Ha\right)}
\global\long\def\MccpH#1{\mathcal{M}_{#1,+}(\Ha)}
\global\long\def\McpHa{\MccpH c}
\global\long\def\oo{\mathbf{0}}
\global\long\def\e{\mathbf{e}_{0}}
\global\long\def\Tt{\widehat{\mathbf{T}}}
\global\long\def\Xx{\mathbf{X}}
\global\long\def\Ttt{\mathcal{T}}
\global\long\def\LHa{\cL\left(\Ha\right)}
\global\long\def\LHax{\cL^{\times}\left(\Ha\right)}
\global\long\def\LHao{\cL\left(\Ha_{1}\right)}
\global\long\def\LMHa{\cL\left(\MHa\right)}
\global\long\def\LMHax{\cL^{\times}\left(\MHa\right)}
\global\long\def\Hao#1{\Ha_{1}\left[#1\right]}
\global\long\def\Haoo#1{\Ha_{1\perp}\left[#1\right]}
\global\long\def\vl#1{\mathcal{V}\left(#1\right)}
\global\long\def\sn{\mathbf{span}\,}
\global\long\def\Xo#1{\Xx_{1}\left[#1\right]}
\global\long\def\Xoo#1{\Xx_{1}^{\perp}\left[#1\right]}
\global\long\def\UHa{\mathfrak{U}\left(\Ha_{1}\right)}
\global\long\def\BoHa{\mathbf{B}_{1}\left(\Ha_{1}\right)}
\global\long\def\Mkf{\mathsf{M}}
\global\long\def\Mkc#1{\Mkf_{c}\left(#1\right)}
\global\long\def\Mkk{\Mkf_{c}}
\global\long\def\bn{\mathbf{n}}
\global\long\def\ba{\mathbf{a}}
\global\long\def\Www#1#2{\mathbf{W}_{#1}\left[#2\right]}
\global\long\def\WsnJ#1{\Www{#1}{s,\bn,J}}
\global\long\def\WlsnJ{\WsnJ{\lambda,c}}
\global\long\def\Wwi#1{\Www{\infty,c}{#1}}
\global\long\def\WinJ{\Wwi{\bn,J}}
\global\long\def\WsnJa#1{\Www{#1}{s,\bn,J;\ba}}
\global\long\def\WlsnJa{\WsnJa{\lambda,c}}
\global\long\def\WLsnJa{\WsnJa{\lambda;\Vf}}
\global\long\def\Uuu#1#2{\mathbf{U}_{#1}\left[#2\right]}
\global\long\def\UsnJ#1{\Uuu{#1}{s,\bn,J}}
\global\long\def\UtsnJ{\UsnJ{\theta,c}}
\global\long\def\UsnJa#1{\Uuu{#1}{s,\bn,J,\ba}}
\global\long\def\UtsnJa{\UsnJa{\theta,c}}
\global\long\def\fio#1{\varphi_{0}\left(#1\right)}
\global\long\def\fiot{\fio{\theta}}
\global\long\def\fion#1{\varphi_{1}\left(#1\right)}
\global\long\def\fiont{\fion{\theta}}
\global\long\def\fioq#1{\varphi_{0}^{2}\left(#1\right)}
\global\long\def\fionq#1{\varphi_{1}^{2}\left(#1\right)}
\global\long\def\sign{\mathrm{sign}\,}
\global\long\def\LT#1{\mathfrak{OT}\left(#1\right)}
\global\long\def\LTp#1{\mathfrak{OT}_{+}\left(#1\right)}
\global\long\def\LTH{\LT{\Ha,c}}
\global\long\def\LTpH{\LTp{\Ha,c}}
\global\long\def\LG#1{\mathfrak{O}\left(#1\right)}
\global\long\def\LGp#1{\mathfrak{O}_{+}\left(#1\right)}
\global\long\def\LGH{\LG{\Ha,c}}
\global\long\def\LGpH{\LGp{\Ha,c}}
\global\long\def\LGm#1{\mathfrak{O}_{-}\left(#1\right)}
\global\long\def\PT#1{\mathfrak{PT}\left(#1\right)}
\global\long\def\PTp#1{\mathfrak{PT}_{+}\left(#1\right)}
\global\long\def\PTH{\PT{\Ha,c}}
\global\long\def\PTpH{\PTp{\Ha,c}}
\global\long\def\PTHC{\PT{\Ha;\Vf}}
\global\long\def\PTpHC{\PTp{\Ha;\Vf}}
\global\long\def\PG#1{\mathfrak{P}\left(#1\right)}
\global\long\def\PGp#1{\mathfrak{P}_{+}\left(#1\right)}
\global\long\def\PGm#1{\mathfrak{P}_{-}\left(#1\right)}
\global\long\def\PGH{\PG{\Ha,c}}
\global\long\def\PGpH{\PGp{\Ha,c}}
\global\long\def\PTmpf#1{\mathfrak{PT}_{\mathsf{fin}}^{\mp}\left(#1\right)}
\global\long\def\PTmpfH{\PTmpf{\Ha,c}}
\global\long\def\PTmp#1{\mathfrak{PT}^{\mp}\left(#1\right)}
\global\long\def\Smxc#1#2{\sqrt{\left|1-\frac{#1^{2}}{#2^{2}}\right|}}
\global\long\def\Sxc#1#2{\sqrt{1-\frac{#1^{2}}{#2^{2}}}}
\global\long\def\smxc#1{\Smxc{#1}c}
\global\long\def\Vf{\mathfrak{V_{f}}}
\global\long\def\Cvp{\widetilde{\Vf}}
\global\long\def\MHac{\mathcal{M}\left(\Ha_{\Vf}\right)}
\global\long\def\Ipm#1#2{\I_{#1}\left[#2\right]}
\global\long\def\PPk{\mathbf{Pk}}
\global\long\def\PkH{\PPk\left(\Ha\right)}
\global\long\def\Haa{\widehat{\Ha}}
\global\long\def\EL#1#2#3{\mathbf{E}_{#1}^{\left[#2,#3\right]}}
\global\long\def\Elsn{\EL{\lambda,c}s{\bn}}
\global\long\def\Elsna{\EL{\lambda,c}s{\bn;\ba}}
\global\long\def\knK{\mathfrak{K}}
\global\long\def\KP#1{\knK\mathfrak{P}\left(#1\right)}
\global\long\def\KPn#1{\knK\mathfrak{P}_{0}\left(#1\right)}
\global\long\def\KPT#1{\knK\mathfrak{P}\mathfrak{T}\left(#1\right)}
\global\long\def\KPTn#1{\knK\mathfrak{P}\mathfrak{T}_{0}\left(#1\right)}
\global\long\def\knU{\mathfrak{U}}
\global\long\def\UP#1{\knU\mathfrak{P}\left(#1\right)}
\global\long\def\UPn#1{\knU\mathfrak{P}_{0}\left(#1\right)}
\global\long\def\UPT#1{\knU\mathfrak{P}\mathfrak{T}\left(#1\right)}
\global\long\def\UPTn#1{\knU\mathfrak{P}\mathfrak{T}_{0}\left(#1\right)}
\global\long\def\UPTmpf#1{\knU\mathfrak{P}\mathfrak{T}_{\mathsf{fin}}^{\mp}\left(#1\right)}
\global\long\def\UPTmp#1{\knU\mathfrak{P}\mathfrak{T}^{\mp}\left(#1\right)}
\global\long\def\xxxv{\widetilde{\x}}
 \global\long\def\trj#1#2{\mathbf{trj}_{#1}\left(#2\right)}
\global\long\def\shp#1#2{#1^{\left\langle +#2\right\rangle }}
\global\long\def\paral#1#2{\parallel_{#2}^{#1}}
\global\long\def\parll{\paral{\,}{\lf}}
\global\long\def\Qkv{\mathrm{(V)}}
\global\long\def\sX{\mathscr{X}}
\global\long\def\sY{\mathscr{Y}}
\global\long\def\sZ{\mathscr{Z}}
\global\long\def\bsf{\boldsymbol{f}}
\global\long\def\bsg{\boldsymbol{g}}
\global\long\def\Qkt{\mathrm{(T)}}
\global\long\def\Qkr{\mathrm{(R)}}
\global\long\def\Co{\mathrm{C_{0}}}
\global\long\def\Qkco{\mathrm{\left(\Co\right)}}
\global\long\def\Con{\mathrm{C_{1}}}
\global\long\def\Qkcon{\mathrm{\left(\Con\right)}}
\global\long\def\Qktv{\mathrm{(TV)}}
\global\long\def\Qkrco{\mathrm{\left(R\Co\right)}}
\global\long\def\CondA{\text{(a)}}
\global\long\def\CondB{\textrm{(b)}}
\global\long\def\CondC{\textrm{(c)}}
\renewcommand{\doctitle}{On Monotonous Separately Continuous Functions}

\title{\doctitle}

\maketitle
%\ifthenelse{\equal{\zbsty}{yes}} 
%  {\Author{\Large Я.І.~Грушка}}
%  {\bigskip\author{\Large Я.І.~Грушка}}
\renewcommand{\docauthor}{Grushka~Ya.I.}
{\author{\large \docauthor}}\medskip{}

\noindent {\large ~ ~ Institute of Mathematics NAS of Ukraine,
Kyiv, created: \textcolor{blue}{January 4, 2018}. }\bigskip{}

\ifthenelse{\equal{\DOI}{}}{}
  {\hfill \textcolor{dgreen}{\large\bf DOI: \href{http://doi.org/\DOI}{\DOI}} ~} \medskip{}

\ifthenelse{\equal{\Preprint}{true}}{
  \ifthenelse{\equal{\DOI}{}}{\newcommand{\PreprDOI}{}}
     {\newcommand{\PreprDOI}{,~ DOI: \DOI}}
  \newcommand{\PreprintSrc}{ArXiv}
  \renewcommand{\docinfoOdd}{\textcolor{dgreen}{\bf Preprint: \PreprintSrc\hfill}}}
  {} % \PreprDOI
\renewcommand{\docinfoEven}{\hfil \textcolor{dgreen}{\bf \doctitle}}
\begin{center}

\newcommand{\AbstrTxt}
{

Let $\TT=\left(\T,\leq\right)$ and $\TT_{1}=\left(\T_{1},\leq_{1}\right)$
be linearly ordered sets and $\sX$ be a topological space. The main
result of the paper is the following: 

If function $\bsf(t,x):\T\times\sX\mapsto\T_{1}$ is continuous in
each variable (``$t$'' and ``$x$'') separately and function $\bsf_{x}(t)=\bsf(t,x)$
is monotonous on $\T$ for every $x\in\sX$, then $\bsf$ is continuous
mapping from $\T\times\sX$ to $\T_{1}$, where $\T$ and $\T_{1}$
are considered as topological spaces under the order topology and
$\T\times\sX$ is considered as topological space under the Tychonoff
topology on the Cartesian product of topological spaces $\T$ and
$\sX$. 

}
{\parbox{15.0cm}{\normalsize \parindent=0.5cm \AbstrTxt}}
\end{center}

\medskip

\markboth{}{}
%\selectlanguage{ukrainian} 
\sloppy
\allowdisplaybreaks 

\large
\selectlanguage{english}

\section{Introduction}

In 1910 W.H. Young had proved the following theorem. 

\begin{theoremA}[see \cite{Young01_monot}] 

Let $f:\R^{2}\mapsto\R$ be separately continuous. If $f(\cdot,y)$
is also monotonous for every $y$, then $f$ is continuous. 

\end{theoremA}

In 1969 this theorem was generalized for the case of separately continuous
function $f:\R^{d}\mapsto\R$ ($d\geq2$):

\begin{theoremB}[see \cite{Krusee01}] 

Let $f:\R^{d+1}\mapsto\R$ ($d\in\N$) be continuous in each variable
separately. Suppose $f\left(t_{1},\cdots,t_{d},\tau\right)$ is monotonous
in each $t_{i}$ separately ($1\leq i\leq d$). Then $f$ is continuous
on $\R^{d+1}$. 

\end{theoremB}

Note that theorems A an B were also mentioned in the overview \cite{Ciesielski_Miller_overview2017}.
In the papers \cite{Mykhaylyuk01_monot,Nesterenko01_monot} authors
investigated functions of kind $\bsf:\T\times\sX\mapsto\R$, where
$\left(\T,\leq\right)$ is linearly ordered set equipped by the order
topology, $\left(\sX,\fS_{\sX}\right)$ is any topological space and
the function $\bsf$ is monotonous relatively the first variable as
well continuous (or quasi-continuous) relatively the second variable.
In particular in \cite{Nesterenko01_monot} it was proven that each
separately quasi-continuous and monotonous relatively the first variable
function $\bsf:\R\times\sX\mapsto\R$ is quasi-continuous relatively
the set of variables. The last result may be considered as the abstract
analog of Young's theorem (Theorem A) for separately quasi-continuous
functions. 

However, we do not know any direct generalization of Theorem A (for
separately continuous and monotonous relatively the first variable
function) in abstract topological spaces at the present time. In the
present paper we prove the generalization of theorems A and B for
the case of (separately continuous and monotonous relatively the first
variable) function $\bsf:\T\times\sX\mapsto\T_{1}$, where $\left(\T,\leq\right)$,
$\left(\T_{1},\leq_{1}\right)$ are linearly ordered sets equipped
by the order topology and $\sX$ is any topological space.

\section{Preliminaries }

Let $\TT=\left(\T,\leq\right)$ be any linearly (ie totally) ordered
set (in the sense of \cite{Birkhoff}). Then we can define the strict
linear order relation on $\T$ such, that for any $t,\tau\in\T$ the
correlation $t<\tau$ holds if and only if $t\leq\tau$ and $t\neq\tau$.
Together with the linearly ordered set $\TT$ we introduce the linearly
ordered set 
\[
\TT_{\pm\infty}=\left(\T\cup\left\{ -\infty,+\infty\right\} ,\leq\right),
\]
where the order relation is extended on the set $\T\cup\left\{ -\infty,+\infty\right\} $
by means of the following clear conventions: 
\begin{description}
\item [{(a)}] $-\infty<+\infty$; 
\item [{(b)}] $\forall t\in\T$ ~$-\infty<t<+\infty$. 
\end{description}
Recall \cite{Birkhoff} that every such linearly ordered set $\TT=\left(\T,\leq\right)$
can be equipped by the natural ``internal'' order topology $\Tpi\left[\TT\right]$,
generated by the base consisting of the open sets of kind:
\begin{equation}
\left(\tau_{1},\tau_{2}\right)=\left\{ t\in\T\,|\:\tau_{1}<t<\tau_{2}\right\} ,\quad\textrm{where}\quad\tau_{1},\tau_{2}\in\T\cup\left\{ -\infty,+\infty\right\} ,\:\tau_{1}<\tau_{2}.\label{eqP:Intervals&Topol_bs}
\end{equation}

Let $\left(\sX,\fS_{\sX}\right)$, $\left(\sY,\fS_{\sY}\right)$ and
$\left(\sZ,\fS_{\sZ}\right)$ be topological spaces, where $\fS_{\mathscr{S}}\subseteq2^{\mathscr{S}}$
is the topology on the topological space $\mathscr{S}$ ($\mathscr{S}\in\left\{ \sX,\sY,\sZ\right\} $).
By $\mathbf{C}(\sX,\sY)$ we denote the collection of all continuous
mappings from $\sX$ to $\sY$. For a mapping $\bsf:\sX\times\sY\mapsto\sZ$
and a point $(x,y)\in\sX\times\sY$ we write 
\[
\bsf^{x}(y):=\bsf_{y}(x):=\bsf(x,y).
\]
Recall \cite{Karlova01} that the mapping $\bsf:\sX\times\sY\mapsto\sZ$
is refereed to as \textbf{\emph{separately continuous}} if and only
if $\bsf^{x}\in\mathbf{C}(\sY,\sZ)$ and $\bsf_{y}\in\mathbf{C}(\sX,\sZ)$
for every point $(x,y)\in\sX\times\sY$ (see also \cite{MasluchenkoVK(NN)2,Mykhaylyuk01_monot,Nesterenko01_monot}).
The set of all separately continuous mappings $\bsf:\sX\times\sY\mapsto\sZ$
is denoted by $\mathbf{CC}\left(\sX\times\sY,\sZ\right)$ \cite{Karlova01,MasluchenkoVK(NN)2,Mykhaylyuk01_monot,Nesterenko01_monot}. 

Let $\TT=\left(\T,\leq\right)$ and $\TT_{1}=\left(\T_{1},\leq_{1}\right)$
be linearly ordered sets. We say that a function $f:\T\mapsto\T_{1}$
is \textbf{\emph{non-decreasing}} (\textbf{\emph{non-increasing}})
on $\T$ if and only if for every $t,\tau\in\T$ the inequality $t\leq\tau$
leads to the inequality $f(t)\leq_{1}f(\tau)$ ($f(\tau)\leq_{1}f(t)$)
correspondingly. Function $f:\T\mapsto\T_{1}$, which is non-decreasing
or non-increasing on $\T$ is called by \textbf{\emph{monotonous}}.

\section{Main Results }

Let $\left(\sX_{1},\fS_{\sX_{1}}\right)$, $\cdots$, $\left(\sX_{d},\fS_{\sX_{d}}\right)$~
($d\in\N$) be topological spaces. Further we consider $\sX_{1}\times\cdots\times\sX_{d}$
as a topological space under the Tychonoff topology $\fS_{\sX_{1}\times\cdots\times\sX_{d}}$
on the Cartesian product of topological spaces $\sX_{1}$, $\cdots$,
$\sX_{d}$. Recall \cite[Chapter 3]{Kelley} that topology $\fS_{\sX_{1}\times\cdots\times\sX_{d}}$
is generated by the base of kind: 
\[
\left\{ U_{1}\times\cdots\times U_{d}\,|\: U_{j}\in\fS_{\sX_{j}}\,\left(\forall j\in\left\{ 1,\cdots,d\right\} \right)\right\} .
\]

\BeginThm  \label{ThmP:YoungGeneralized}

Let $\TT=\left(\T,\leq\right)$ and $\TT_{1}=\left(\T_{1},\leq_{1}\right)$
be linearly ordered sets and $\left(\sX,\fS_{\sX}\right)$ be a topological
space. 

If $\bsf\in\mathbf{CC}\left(\T\times\sX,\T_{1}\right)$ and function
$\bsf_{x}(t)=\bsf(t,x)$ is monotonous on $\T$ for every $x\in\sX$,
then $\bsf$ is continuous mapping from the topological space $\left(\T\times\sX,\fS_{\,\T\times\sX}\right)$
to the topological space $\left(\T_{1},\Tpi\left[\TT_{1}\right]\right)$.

\EndThm  

\BeginProof  

Consider any ordered pair $\left(t_{0},x_{0}\right)\in\T\times\sX$.
Take any open set $V\subseteq\T_{1}$ such, that $\bsf\left(t_{0},x_{0}\right)\in V$.
Since the sets of kind (\ref{eqP:Intervals&Topol_bs}) form the base
of topology $\Tpi\left[\TT_{1}\right]$, there exist $\tau_{1},\tau_{2}\in\T_{1}\cup\left\{ -\infty,+\infty\right\} $
such, that $\tau_{1}<_{1}\bsf\left(t_{0},x_{0}\right)<_{1}\tau_{2}$
and $\left(\tau_{1},\tau_{2}\right)\subseteq V$, where $<_{1}$ is
the strict linear order, generated by (non-strict) order $\leq_{1}$
(on $\T_{1}\cup\left\{ -\infty,+\infty\right\} $). The function $\bsf$
is separately continuous. So, since the sets of kind (\ref{eqP:Intervals&Topol_bs})
form the base of topology $\Tpi\left[\TT\right],$ there exist $t_{1},t_{2}\in\T\cup\left\{ -\infty,+\infty\right\} $
such, that 
\begin{align}
 & t_{1}<t_{0}<t_{2}\quad\textrm{and}\quad\label{eqP:(t1,t2):props1}\\
 & \bsf\left(\left(t_{1},t_{2}\right),x_{0}\right)\subseteq\left(\tau_{1},\tau_{2}\right).\label{eqP:(t1,t2):props2}
\end{align}

Further we need the some additional denotations. 
\begin{itemize}
\item In the case, where $\left(t_{1},t_{0}\right)\neq\emptyset$ we choose
any element $\alpha_{1}\in\T$ such that $t_{1}<\alpha_{1}<t_{0}$
and denote $\widetilde{\alpha}_{1}:=\alpha_{1}$. In the opposite
case we denote $\alpha_{1}:=t_{0}$, $\widetilde{\alpha}_{1}:=t_{1}$. 
\item In the case $\left(t_{0},t_{2}\right)\neq\emptyset$ we choose any
element $\alpha_{2}\in\T$ such that $t_{0}<\alpha_{2}<t_{2}$ and
denote $\widetilde{\alpha}_{2}:=\alpha_{2}$. In the opposite case
we denote $\alpha_{2}:=t_{0}$, $\widetilde{\alpha}_{2}:=t_{2}$. 
\end{itemize}
It is not hard to verify, that in the all cases the following conditions
are performed: 
\begin{align}
 & \alpha_{1},\alpha_{2}\in\T,\quad\widetilde{\alpha}_{1},\widetilde{\alpha}_{2}\in\T\cup\left\{ -\infty,+\infty\right\} ;\nonumber \\
 & \alpha_{1}\leq\alpha_{2};\nonumber \\
 & \widetilde{\alpha}_{1}<\widetilde{\alpha}_{2};\nonumber \\
 & \left[\alpha_{1},\alpha_{2}\right]\subseteq\left(t_{1},t_{2}\right),\quad\textrm{where}\;\left[\alpha_{1},\alpha_{2}\right]=\left\{ t\in\T\:|\:\alpha_{1}\leq t\leq\alpha_{2}\right\} ;\label{eqP:[al1,al2]subs(t1,t2)}\\
 & t_{0}\in\left(\widetilde{\alpha}_{1},\widetilde{\alpha}_{2}\right)\subseteq\left[\alpha_{1},\alpha_{2}\right].\label{eqP:t0(al1x,al2x)[al1,al2]}
\end{align}

According to (\ref{eqP:[al1,al2]subs(t1,t2)}), $\alpha_{1},\alpha_{2}\in\left(t_{1},t_{2}\right)$.
Hence, according to (\ref{eqP:(t1,t2):props2}), interval $\left(\tau_{1},\tau_{2}\right)$
is the open neighborhood of the both points $\bsf\left(\alpha_{1},x_{0}\right)$
and $\bsf\left(\alpha_{2},x_{0}\right)$. Since the function $\bsf$
is separately continuous on $\T\times\sX$, then there exist an open
neighborhood $U\in\fS_{\sX}$ of the point $x_{0}$ (in the space
$\sX$) such, that:
\begin{align}
\bsf\left(\alpha_{1},U\right) & \subseteq\left(\tau_{1},\tau_{2}\right);\label{eqP:f(al1,U)subs(tau12)}\\
\bsf\left(\alpha_{2},U\right) & \subseteq\left(\tau_{1},\tau_{2}\right).\label{eqP:f(al2,U)subs(tau12)}
\end{align}
 The set $\left(\widetilde{\alpha}_{1},\widetilde{\alpha}_{2}\right)\times U$
is the open neighborhood of the point $\left(t_{0},x_{0}\right)$
in the topology $\fS_{\,\T\times\sX}$ of the space $\T\times\sX$.
Now our aim is to prove that 
\begin{equation}
\forall\,\left(t,x\right)\in\left(\widetilde{\alpha}_{1},\widetilde{\alpha}_{2}\right)\times U\;\left(\bsf\left(t,x\right)\in\left(\tau_{1},\tau_{2}\right)\subseteq V\right).\label{eqP:forAall(t,x)f(t,x)in(tau12)}
\end{equation}
 So, chose any point $\left(t,x\right)\in\left(\widetilde{\alpha}_{1},\widetilde{\alpha}_{2}\right)\times U$.
According to the condition (\ref{eqP:t0(al1x,al2x)[al1,al2]}), we
have $\left(t,x\right)\in\left[\alpha_{1},\alpha_{2}\right]\times U$,
that is $\alpha_{1}\leq t\leq\alpha_{2}$ and $x\in U$. In accordance
with (\ref{eqP:f(al1,U)subs(tau12)}), (\ref{eqP:f(al2,U)subs(tau12)}),
we have $\bsf\left(\alpha_{1},x\right)\in\left(\tau_{1},\tau_{2}\right)$
and $\bsf\left(\alpha_{2},x\right)\in\left(\tau_{1},\tau_{2}\right)$.
Hence, since the function $\bsf_{x}(\cdot)=\bsf(\cdot,x)$ is monotonous
on $\T$ and $\alpha_{1}\leq t\leq\alpha_{2}$, we deduce $\bsf\left(t,x\right)\in\left(\tau_{1},\tau_{2}\right)\subseteq V$.
Thus, the correlation (\ref{eqP:forAall(t,x)f(t,x)in(tau12)}) is
proven. Hence, the function $\bsf$ is continuous in (every) point
$\left(t_{0},x_{0}\right)\in\T\times\sX$.          \EndProof  

Theorem~A is a consequence of Theorem~\ref{ThmP:YoungGeneralized}
in the case $\T=\sX=\R$, where $\R$ is considered together with
the usual linear order on the field of real numbers and usual topology. 

\BeginNasl \label{NaslP:YoungGeneralized}

Let $\TT_{0}=\left(\T_{0},\leq_{0}\right)$, $\TT_{1}=\left(\T_{1},\leq_{1}\right)$,
$\cdots$, $\TT_{d}=\left(\T_{d},\leq_{d}\right)$ ($d\in\N$) be
linearly ordered sets, and $\left(\sX,\fS_{\sX}\right)$ be a topological
space. 

If the function $\bsf:\T_{1}\times\cdots\times\T_{d}\times\sX\mapsto\T_{0}$
is continuous in each variable separately and $f\left(t_{1},\cdots,t_{d},\tau\right)$
is monotonous in each $t_{i}$ separately ($1\leq i\leq d$) then
$\bsf$ is a continuous mapping from the topological space $\left(\T_{1}\times\cdots\times\T_{d}\times\sX,\,\fS_{\,\T_{1}\times\cdots\times\T_{d}\times\sX}\right)$
to the topological space $\left(\T_{0},\Tpi\left[\TT_{0}\right]\right)$. 

\EndNasl  

\BeginProof  

We will prove this corollary by induction. For $d=1$ the corollary
is true by Theorem~\ref{ThmP:YoungGeneralized}. Assume, that the
corollary is true for the number $d-1$, where $d\in\N$, $d\geq2$.
Suppose, that function $\bsf:\T_{1}\times\cdots\times\T_{d}\times\sX\mapsto\T_{0}$
is satisfying the conditions of the corollary. Then we may consider
this function as a mapping from $\T_{1}\times\sX_{(d)}$ to $\T_{0}$,
where $\sX_{(d)}=\T_{2}\times\cdots\times\T_{d}\times\sX$. According
to inductive hypothesis, function $\bsf\left(t_{1},\cdot\right)$
is continuous on $\sX_{(d)}$ for every fixed $t_{1}\in\T_{1}$. So
$\bsf$ is a separately continuous mapping from $\T_{1}\times\sX_{(d)}$
to $\T_{0}$. Moreover, $\bsf$ is monotonous relatively the first
variable (by conditions of the corollary). Hence, by Theorem \ref{ThmP:YoungGeneralized},
$\bsf$ is continuous on $\T_{1}\times\sX_{(d)}$.        \EndProof  

Theorem~B is a consequence of Corollary~\ref{NaslP:YoungGeneralized}
in the case $\T_{0}=\T_{1}=\cdots=\T_{d}=\sX=\R$, where $\R$ is
considered together with the usual linear order on the field of real
numbers and usual topology. In the case $\T_{0}=\R$, $\T_{j}=\left(a_{j},b_{j}\right)$,
$\sX=\left(a_{d+1},b_{d+1}\right)$ where $a_{j},b_{j}\in\R$ and
$a_{j}<b_{j}$ ($j\in\left\{ 1,\cdots,d+1\right\} $) and intervals
$\left(a_{j},b_{j}\right)$ are considered together with the usual
linear order and topology, induced from the field of real numbers,
we obtain the following corollary. 

\BeginNasl  \label{NaslP:YoungGeneralized1}

If the function $f:\left(a_{1},b_{1}\right)\times\cdots\times\left(a_{d},b_{d}\right)\times\left(a_{d+1},b_{d+1}\right)\mapsto\R$
($d\in\N$) is continuous in each variable separately and $f\left(t_{1},\cdots,t_{d},\tau\right)$
is monotonous in each $t_{i}$ separately ($1\leq i\leq d$) then
$f$ is a continuous mapping from $\left(a_{1},b_{1}\right)\times\cdots\times\left(a_{d+1},b_{d+1}\right)$
to $\R$. 

\EndNasl  

\BeginRmk 

In fact in the paper \cite{Krusee01} the more general result was
formulated, in comparison with Theorem~B. Namely the author of \cite{Krusee01}
had considered the real valued function $f\left(t_{1},\cdots,t_{d},\tau\right)$
defined on an open set $G\subseteq\R^{d+1}$, $d\in\N$ such, that
$f$ is continuous in each variable separately and monotonous in each
$t_{i}$ separately ($1\leq i\leq d$). But this result of \cite{Krusee01}
can be delivered from Corollary~\ref{NaslP:YoungGeneralized1}, because
for each point $\mathbf{t}=\left(t_{1},\cdots,t_{d},\tau\right)\in G$
in the open set $G$ there exists the set of intervals $\left(a_{1},b_{1}\right),\cdots,\left(a_{d+1},b_{d+1}\right)$
such, that $\mathbf{t}\in\left(a_{1},b_{1}\right)\times\cdots\times\left(a_{d+1},b_{d+1}\right)\subseteq G$. 

\EndRmk  

\ifthenelse{\equal{\Preprint}{true}} {

\selectlanguage{english}%

\paragraph*{Notes on applications in abstract kinematics. }

Main results of the paper are expected to be applied in the theory
of universal kinematics for establishing some additional properties
of coordinate transform operators (between reference frames), separately
continuous relatively space and time variables. For more details about
the theory of universal kinematics see \cite{MyTmm_2017_0(v2),MyTmm10}
and other papers, reference to which you can find in \cite{MyTmm_2017_0(v2),MyTmm10}. \selectlanguage{american}%

}{}

\providecommand{\bysame}{\leavevmode\hbox to3em{\hrulefill}\thinspace}


\begin{thebibliography}{10}

\bibitem{Birkhoff}
G.~Birkhoff, \emph{{Lattice theory}}, {Third edition. American Mathematical
  Society Colloquium Publications, Vol. XXV}, American Mathematical Society,
  Providence, R.I., New York, 1967.

\bibitem{Ciesielski_Miller_overview2017}
K.C. Ciesielski and D.~Miller, \emph{{A Continuous Tale on Continuous and
  Separately Continuous Functions}}, Real Analysis Exchange \textbf{41} (2016),
  no.~1, 19--54,
  \url{http://www.jstor.org/stable/10.14321/realanalexch.41.1.0019}.

\bibitem{MyTmm10}
Y.~Grushka, \emph{{Kinematic changeable sets with given universal coordinate
  transforms}}, Zb. Pr. Inst. Mat. NAN Ukr. \textbf{12} (2015), no.~1, 74--118.

\bibitem{MyTmm_2017_0(v2)}
Y.~Grushka, \emph{{Draft introduction to abstract kinematics. (Version 2.0)}},
  Preprint: ResearchGate, 2017, pp.~1--208,
  \url{https://doi.org/10.13140/RG.2.2.28964.27521}.

\bibitem{Karlova01}
O.~Karlova, V.~Mykhaylyuk, and O.~Sobchuk, \emph{{Diagonals of separately
  continuous functions and their analogs}}, Topology and its Applications
  \textbf{160} (2013), no.~1, 1--8,
  \href{http://dx.doi.org/10.1016/j.topol.2012.09.003}{\path{doi:10.1016/j.topol.2012.09.003}}.

\bibitem{Kelley}
J.L. Kelley, \emph{{General Topology}}, { University series in higher
  mathematics}, Van Nostrand, 1955.

\bibitem{Krusee01}
R.L. Krusee and J.J. Deely, \emph{{Joint Continuity of Monotonic Functions}},
  The American Mathematical Monthly \textbf{76} (1969), no.~1, 74--76,
  \url{http://www.jstor.org/stable/2316804}.

\bibitem{Mykhaylyuk01_monot}
V.~{Mykhajlyuk}, \emph{{The Baire classification of separately continuous and
  monotone functions.}}, {Scientific Herald of Yuriy Fedkovych Chernivtsi
  National University} \textbf{349} (2007), 95--97 (Ukrainian).

\bibitem{Nesterenko01_monot}
V.~Nesterenko, \emph{{Joint properties of functions which monotony with respect
  to the first variable.}}, {Mathematical Bulletin of Taras Shevchenko
  Scientific Society} \textbf{6} (2009), 195--201 (Ukrainian).

\bibitem{MasluchenkoVK(NN)2}
H.~Voloshyn, V.~Maslyuchenko, and O.~Maslyuchenko, \emph{{On layer-wise uniform
  approximation of separately continuous functioins by polynomials}},
  Mathematical Bulletin of Taras Shevchenko Scientific Society \textbf{10}
  (2013), 135--158 (Ukrainian).

\bibitem{Young01_monot}
W.~Young, \emph{{A note on monotone functions}}, The Quarterly Journal of Pure
  and Applied Mathematics (Oxford Ser.) \textbf{41} (1910), 79--87.

\end{thebibliography}
\end{document}